\documentclass[a4paper,10pt]{article}
\usepackage{amsfonts}
\usepackage{amstext}
\usepackage{amssymb,amsmath}
\usepackage{graphicx}
\usepackage{epstopdf}
\usepackage[T1]{fontenc}
\usepackage[utf8]{inputenc}
\usepackage{authblk}
\usepackage{bm}

\newtheorem{theorem}{Theorem}
\newtheorem{lemma}{Lemma}

\newtheorem{proof}{Proof}
\newtheorem{proposition}{Proposition}

\def\R{\Bbb R}
\def\E{\Bbb E}\def\P{\Bbb P}
\def\0{\bold 0}
\def\1{\bold 1}

\def\a{\bold a}
\def\b{\bold b}

\def\t{\bold t}
\def\x{\bold x}

\def\y{\bold y}

\def\U{\bold U}
\def\V{\bold V}

\def\J{\bold J}

\def\RR{\bold R}
\def\W{\bold W}

\def\CC{\bold C}

\def\thetav{\bm{\theta}}

\def\part{\cal P}

\title{Estimating parameters of a directed weighted graph model
with beta-distributed edge-weights}

\author[1]{Marianna Bolla \thanks{marib@math.bme.hu}}

\author[1]{Ahmed Elbanna \thanks{ahmed@math.bme.hu}}

\author[1,2]{J\'ozsef Mala \thanks{jmala@math.bme.hu}}

\affil[1]{Institute of Mathematics, Budapest University of Technology and
Economics \\
1111. Budapest, M\H uegyetem rkp. 3, Hungary}

\affil[2]{ELTE E\"otv\"os Lor\'and University, Institute of Mathematics, Budapest, Hungary}
\begin{document}

\maketitle

\section*{Abstract}

We introduce a directed, weighted random graph model, where the edge-weights
are independent and beta-distributed with parameters depending on their
endpoints. We will show that the row- and column-sums of the transformed
edge-weight matrix are sufficient statistics for the parameters, and
use the theory of exponential families to prove that the ML estimate of
the parameters exists and is unique. Then an algorithm to find this
estimate is introduced together with convergence proof that uses properties
of the digamma function. Simulation results and applications are also presented.

\noindent
\textbf{Keywords}: exponential family,
sufficient statistics, ML estimation, digamma function, successive
approximation

\noindent
\textit{MSC2010}: 62F10, 62B05. 

\section{Introduction}\label{intro}

The theory of ML estimation in the following types of
exponential family random graph models
has frequently been investigated in the last decade,
see, e.g.,~\cite{Bolla,Chatterjee,Hillar13,Yan,Wainwright}.
The graph has $n$ vertices, and the adjacency relations between them are
given by the $n\times n$ random edge-weight matrix $\W =(w_{ij} )$
of zero diagonal.
If $\W$ is symmetric, then we have an  undirected
graph; otherwise, our graph is directed, where $w_{ij}$
is the nonnegative weight assigned to the $i\to j$ edge according to the
model.
We assume that the edge-weights (above or out of the main diagonal) are
completely independent (but their distribution usually depends on
different parameters),
and have an exponential family distribution $\P_{\thetav}$.
So the likelihood function has the general form
\begin{equation}\label{ef}
  L_{\thetav } (\W )
= e^{<\thetav , \t (\W ) > \, - \, Z (\thetav ) } \cdot h(\W ) ,
\end{equation}
with the \textit{canonical parameter} $\thetav$,
\textit{log-partition (cumulant) function} $Z (\thetav )$, and
\textit{canonical sufficient statistic} $\t$.
In these random graph models, components of $\t =\t (\W )$ are the row-sums
and/or column-sums of $\W$ or some $\W$-related matrix, i.e., they are
vertex-degrees or in- and out-degrees of the observed undirected or directed,
weighted or unweighted graph (in the weighted case, the edge-weights
may undergo a suitable transformation).
Also, $h (\W )$ is usually 1 over the support of the
likelihood function, indicating that given the canonical sufficient statistics,
the joint distribution of the entries
is uniform (microcanonical) in these models.

To make inferences on the parameters, typically we have only one observation
for the graph.
It may seem that it is a one-element sample, but there are the adjacencies
that form the sample; the number of them is $n\choose 2$
in the undirected, and $n(n-1)$ in the directed case. The number of parameters,
contained in $\thetav$, is $n$ in the undirected and $2n$ in the directed case.
The parameters
can be considered as affinities or potentials of the vertices to make
ties in the undirected, and to emanate or adsorb edges in the
directed case. It is important that we divide the components of the
canonical parameter $\thetav$ of the underlying distribution
of the $ij$ or $i\to j$ edge between the connected vertices,
like $\alpha_i +\alpha_j$ in the undirected and
$\alpha_i +\beta_j$ in the directed case ($i\ne j$),
see~\cite{Bolla,Chatterjee,Yan}.

In \textit{regular exponential families} ($\Theta$ is open), the ML equation
$
 \nabla_{\thetav} \ln L_{\thetav} (\W ) = \0
$
is equivalent to
\begin{equation}\label{Z}
 \nabla_{\thetav} Z(\thetav ) = \t  .
\end{equation}
Since $\nabla_{\thetav} Z(\thetav ) = \E_{\thetav } \t$,
the ML equation~(\ref{Z}) means that the canonical sufficient statistic
is made equal
to its expectation. But when is it possible? Now we briefly summarize
existing theoretical  results on this issue.
Let ${\cal M} = \{ \E_{\thetav } \t \, : \, \thetav \in \Theta \}$ denote the
so-called \textit{mean parameter space} in the model; it is necessarily convex.
Let  ${\cal M}^0$ denote its interior.
When the canonical statistic is also complete,
and hence, minimal sufficient, the representation~(\ref{ef}) is
\textit{minimal} (i.e., the model is not overparametrized).

\begin{proposition}[Proposition 3.2 of~\cite{Wainwright}]
In exponential family, the gradient mapping $\nabla Z: \, \Theta \to {\cal M}$
is one-to-one if and only if the exponential family representation is
minimal.
\end{proposition}\label{pr1}
\begin{proposition}[Theorem 3.3 of~\cite{Wainwright}]\label{pr2}
In a minimal exponential family, the gradient mapping $\nabla Z$ is onto
${\cal M}^0$.
\end{proposition}
By Propositions~\ref{pr1} and~\ref{pr2}, any parameter in ${\cal M}^0$
is uniquely realized by the $\P_{\theta }$ distribution
for some $\theta \in \Theta$. Also, in a regular and minimal
exponential family, ${\cal M}$ is an open set and  is identical
to ${\cal M}^0$.

As the ML estimate of $\thetav$  is the solution of~(\ref{Z}),
we have the following.

\begin{proposition}[Proposition 5 of~\cite{Yan}]\label{pr3}
Assume, the (canonical) parameter space $\Theta$ is open. Then there exists a
solution ${\hat \thetav} \in \Theta$ to the ML equation
$\nabla_{\thetav} Z(\thetav ) = \t$ if and only if $\t \in {\cal M}^0$;
further, if such a solution exists, it is also unique.
\end{proposition}

Note that in regular and minimal exponential families, ${\cal M}^0$ is also the
interior of ${\cal T}$, which is the convex hull of all possible values of
$\t$, see, e.g.,~\cite{Chatterjee,Lauritzen}.
In the case of discrete distributions, it frequently happens that the boundary
of $\cal T$ has positive measure. For instance, the so-called threshold graphs
are located on the boundary of the polyhedron, determined by the
Erd\H os--Gallai conditions, in the model
of~\cite{Chatterjee} which uses Bernoulli distributed entries.
However, in the case of an
absolutely continuous $\P_{\thetav}$ distribution,
the boundary of $\cal T$
has zero Lebesgue measure, and so, probability zero with respect to the
 $\P_{\thetav}$ measure.
Therefore, in view of Proposition~\ref{pr3},
the ML equation has a unique solution with probability 1.

The organization of the paper is as follows. In Section~\ref{model}, we
introduce a model for directed edge-weighted graphs and prove that
a unique ML estimate of the parameters exists.
In Section~\ref{algorithm}, we define an iterative algorithm to find
this solution, and prove its convergence with a convenient starting.
In Section~\ref{application},
the algorithm is applied to randomly generated and real-word data.
In Appendix A, properties of the digamma function, whereas in Appendix B,
the boundary of our
$\cal M$ is discussed.  The long proof of the main
convergence theorem of the iteration algorithm, introduced in
Section~\ref{algorithm}, is presented in Appendix C.

We remark that edge-weighted graphs of uniformly bounded edge-weights are
prototypes of real-world networks, see e.g.,~\cite{Bolla13}.
Without loss of generality, if the
edge-weights are transformed into the [0,1] interval, the
\textit{beta-distribution}
for them, with varying parameters, is capable to model a wide range of
possible probability densities on them. This indicates the soundness of
the model to be introduced in Section~\ref{model}.

\section{A  random graph model with beta-distributed
edge-weights}\label{model}

Let $\W =(w_{ij})$ be the $n\times n$ (usually not symmetric)
edge-weight matrix of a random directed graph
on $n$ vertices: $w_{ii} =0$ $(i=1,\dots ,n)$ and $w_{ij} \in [0,1]$ is the
weight of the $i\to j$ edge $(i\ne j)$. Our model is the following:
the $i\ne j$ weight obeys a beta-distribution  with parameters $a_i >0$ and
$b_j >0$. The parameters are collected in $\a =(a_1 ,\dots a_n )$ and
$\b= (b_1 ,\dots b_n )$, or briefly, in $\thetav = (\a ,\b )$.
Here $a_i$ can be thought of as the potential of the vertex $i$ to send
messages out, and $b_i$ is its resistance to receive messages in.

The likelihood function is factorized as
$$
\begin{aligned}
L_{\a ,\b } (\W ) &=
 \prod_{i\ne j} \frac{\Gamma (a_i  +b_j )}{\Gamma (a_i ) \Gamma (b_j )}
 w_{ij}^{a_i -1} (1-w_{ij})^{b_j -1} \\
& = C(\a ,\b )\prod_{i\ne j} \exp
 [ (a_i -1 ) \ln w_{ij} + (b_j -1) \ln (1-w_{ij}) ] \\
& = \exp
 [ \sum_{i=1}^n (a_i -1 ) \sum_{j\ne i} \ln w_{ij} +
   \sum_{j=1}^n (b_j -1 ) \sum_{i\ne j} \ln (1-w_{ij} ) -Z (\a , \b)],
\end{aligned}
$$
where $C(\a ,\b )$ is the normalizing constant, and
$Z (\a , \b) = -\ln C(\a ,\b )$ is the log-partition (cumulant) function.
Since the likelihood function depends on $\W$ only through the row-sums of
the $n\times n$ matrix $\U =\U (\W )$ of general entry $\ln w_{ij}$ and the
column-sums  of the $n\times n$ matrix $\V = \V (\W)$
of general entry $\ln (1-w_{ij} )$,
by the Neyman--Fisher factorization theorem, the row-sums
$R_1 ,\dots ,R_n$ of $\U$ and column-sums  $C_1 ,\dots ,C_n$ of $\V$ are
sufficient statistics for the parameters. Moreover,
$\t =(\RR ,\CC ) = (R_1 ,\dots ,R_n , C_1 ,\dots ,C_n)$ is the canonical
sufficient statistic, which is also minimal.
Note that $\U$ contains the log-weights of the original graph,
while $\V$ contains the
the log-weights of the complement graph of edge-weight matrix $\overline \W $
with entries $1-w_{ij}$ $(i\ne j )$.
The first factor in the Neyman--Fisher factorization
(including gamma-functions)
depends only on the parameters and on the sample through these
sufficient statistics,
whereas the seemingly not present other
factor -- which would merely depend on $\W$ -- is constantly
1, indicating that the conditional
joint distribution of the entries, given the row- and column-sums of the
log-weight and log-complement matrix
is uniform (microcanonical) in this model.
So under the conditions on the margins of $\U$ and $\V$, the
directed graphs coming from the above model are uniformly distributed.

The system of likelihood equations is obtained by making the derivatives of
$L_{\a ,\b } (\W )$ with respect to the parameters equal to 0:
\begin{equation}\label{le}
\begin{aligned}
\frac{\partial  L_{\a ,\b } (\W )}{\partial a_i } &=
 \sum_{j\ne i}  \psi (a_i  +b_j ) - (n-1) \psi (a_i) +R_i =0,
 \quad i=1,\dots ,n; \\
\frac{\partial L_{\a ,\b } (\W )}{\partial b_j } &=
 \sum_{i\ne j}  \psi (a_i  +b_j ) - (n-1) \psi (b_j) +C_j =0 ,
\quad j=1,\dots ,n .
\end{aligned}
\end{equation}
Here $\psi (x)=\frac{\partial \ln \Gamma (x)}{\partial x} =
\frac{\Gamma' (x)}{\Gamma (x)}$ for $x>0$
is the \textit{digamma function}. For its properties, see Appendix A.

To apply the theory of Section~\ref{intro},
we utilize that the parameter space
$\Theta \subset \R_{+}^{2n}$ is open, akin to the canonical
parameter space, $(-1 ,\infty )^{2n}$.
Note that the canonical parameter  is, in fact,
$(\a' ,\b' )={\thetav}' =\thetav -\1$,
where $\1\in\R^{2n}$ is the vector of all 1 coordinates.
With it, the log-partition function  is
$$
 Z ({\a}', {\b}' )=
 -\sum_{j\ne i} \Gamma ({ a_i}'  +{ b_j}' +2 ) +
  \sum_{j\ne i} \Gamma ({a_i}' +1)
  +\sum_{i\ne j} \Gamma ({b_j}' +1) .
$$
In view of~(\ref{Z}), the ML equation is equivalent to
$$
\begin{aligned}
 \frac{\partial Z ({\a}', {\b}' )}{\partial {a_i}' } &=
 -\sum_{j\ne i} \psi ({a_i}'  +{b_j}' +2 ) +
  (n-1) \psi ({a_i}' +1)  =R_i , \quad i=1,\dots ,n; \\
 \frac{\partial Z ({\a}', {\b}' )}{\partial {b_j}' } &=
 -\sum_{i\ne j} \psi ({a_i}'  +{b_j }' +2 ) +
  (n-1) \psi ({b_j}' +1)  =C_j  , \quad i=1,\dots ,n  .
\end{aligned}
$$
But this system of equations is the same as~(\ref{le}), in terms of
the parameter $\thetav'$ instead of $\thetav$.

In view of  Section~\ref{intro},
the mean parameter space ${\cal M}$ consists of parameters
$(A_1 ,\dots A_n ,B_1 ,\dots ,B_n)$ obtained by the gradient mapping, that is,
\begin{equation}\label{mean}
\begin{aligned}
 A_i = A_i (\a ,\b ) =- \sum_{j\ne i}  \left[ \psi ({a}_i  +{b}_j ) -
  \psi ({a}_i  ) \right] , \quad i=1,\dots ,n; \\
 B_j = B_j (\a ,\b ) =- \sum_{i\ne j} \left[ \psi ({a}_i  +{b}_j  ) -
   \psi ({b}_j ) \right]  , \quad j=1,\dots ,n .
\end{aligned}
\end{equation}
$\cal M$ is an open set, whose boundary
is determined by the limit properties between the digamma and the log functions,
see Appendix B for details. There we also find a correspondence between
the points on the boundary of $\cal M$ and those on the boundary
of the  convex hull $\cal T$
of the possible sufficient statistics $\t =(\RR ,\CC )$ within $\R_{-}^{2n}$.
It is interesting that while the boundary points of $\cal M$ do not belong
to the open set $\cal M$, the boundary points of $\cal T$ do belong to $\cal T$,
and can be realized as row- and column-sums of the $\U (\W )$ and $\V (\W )$
matrices with a $\W$ of off-diagonal entries in  $(0,1)$.
However, this boundary has 0 probability, and so, any canonical
sufficient statistic $\t$ of the observed graph  is in $\cal M$,
with probability 1.
Therefore, by Proposition~\ref{pr3},  we can state the following.
\begin{theorem}\label{t1}
The system of the ML equations~(\ref{le}) has a unique
solution $\hat \thetav =({\hat \a}, {\hat \b})$, with probability 1.
\end{theorem}

Later we will use the following
trivial upper bound for the sum of row- and column-sums (of the $\U$ and $\V$
 matrices):
\begin{equation}\label{joska1}
 \sum_{i=1}^n R_i + \sum_{j=1}^n C_j =
 \sum_{i=1}^n \sum_{j\ne i} \ln w_{ij}
 + \sum_{j=1}^n \sum_{i\ne j} \ln (1-w_{ij}) =
  \sum_{i\ne j} \ln [w_{ij} (1-w_{ij} )] \le -2\ln 2 \, n (n-1)
\end{equation}
due to rearranging the terms and the relation $w_{ij} (1-w_{ij} )\leq 1/4$ for
$w_{ij} \in [0,1]$ with equality if and only if $w_{ij}=\frac12$ $(i\ne j)$.
For finer estimates see Appendix B.

Also note that the Hessian of the system of ML equations (consisting of
the second order partial derivatives of $L_{\thetav}$) at
$\hat \thetav$) does not contain the  sufficient statistics any more, therefore
the negative of it is the Fisher-information matrix at $\hat \thetav$.
Because of
the regularity conditions, the information matrix is positive, and so,
the Hessian is negative definite.
This is also an indication of the existence of a unique ML estimate.

\section{Iteration algorithm to find the parameters}\label{algorithm}

To use a fixed point iteration,
now we rewrite the system of likelihood equations in the form
$\thetav =f(\thetav )$, where $\thetav = (\a ,\b )$, as follows:
\begin{equation}\label{fixes}
\begin{aligned}
 a_i &= {\psi}^{-1} \left[ \frac1{n-1} R_i +\frac1{n-1}
 \sum_{j\ne i} \psi  (a_i  +b_j ) \right]
  =:g_i (\a ,\b ), \quad i=1,\dots ,n \\
 b_j &= {\psi}^{-1} \left[ \frac1{n-1} C_j +\frac1{n-1}
 \sum_{i\ne j} \psi  (a_i  +b_j ) \right]
  =:h_j (\a ,\b ), \quad j=1,\dots ,n .
\end{aligned}
\end{equation}
Here $g_i$'s and $h_j$'s are the coordinate functions of
$f =(g,h) :\R^{2n} \to \R^{2n}$.
Then, starting at $\thetav^{(0)}$, we use the successive approximation
$\thetav^{(it )} :=f (\thetav^{(it -1)}$ for $it =1,2,\dots$, until convergence.
Now the the statement of convergence of the above iteration to the theoretically
guaranteed unique $\hat \thetav$ (see Theorem~\ref{t1}) follows.

\begin{theorem}\label{t2}
Let ${\hat \thetav} =({\hat \a },{\hat \b} )$ be the unique solution of the
ML equation~(\ref{le}). Then the above
mapping $f=(g,h)$ is a contraction in some closed
neighborhood $K$ of $\hat \thetav$, and so, starting at any
$\thetav^{(0)} \in K$, the fixed point of the iteration
$\thetav^{(it )} =f (\thetav^{(it -1)})$ exists and is $\hat \thetav$.
\end{theorem}
The prof of this theorem is to be found in Appendix C.

Since $K$ is only theoretically guaranteed, we need some practical
considerations about the choice of $\thetav^{(0 )}$, which should be
adapted to the sufficient statistics.
In the sequel, for two vectors $\x =(x_1,\dots,x_n)$,
$\y =(y_1,\dots,y_n)$  we use the notation
$\x >\y$ if $x_i>y_i$ for each $i=1,\dots,n$. Likewise,
$\x \ge \y$ is the shorthand for $x_i\ge y_i$ for each $i=1,\dots,n$.

Recall that $f=(g,h)$ is the mapping (\ref{fixes}) of the fixed point
iteration, and
${\hat\thetav}=({\hat\a},{\hat\b})>{\mathbf 0}$ is the (only) solution of the
equation
$f(\thetav)=\thetav$, where $\0\in\R^{2n}$ is the vector of all 0 coordinates.

\begin{proposition}\label{joskaprop2}
Let
\begin{equation}\label{M}
M:=\max\left\{ \max_{i\in \{ 1,\dots,n \} }\left(-\frac{R_i}{n-1}\right),
\max_{i\in \{ 1,\dots,n \} }\left(-\frac{C_i}{n-1}\right)\right\}
\end{equation}
and $\varepsilon>0$ be the (only) solution of the equation
$\psi(2x)-\psi(x)=M$. Then
$({\hat\a},{\hat\b})\geq\varepsilon\1$.
\end{proposition}

\noindent \textbf{Proof.}
In view of  (\ref{joska1}) we have that $M\geq\ln 2$. Since
equality in (\ref{joska1}) is attained with probability 0, we have that
$M>\ln 2$ with probability 1. Therefore, by
Lemma~\ref{joskaprop1} of Appendix A,
there exists an $\varepsilon$, with probability 1, such that
$\psi(2\varepsilon )-\psi(\varepsilon )=M$.

Without loss of generality we can assume that
$$
{\hat a_{i_0}}=\min\left\{ \min_{i\in \{ 1,\dots,n \} }{\hat a_i},
 \min_{i\in \{ 1,\dots,n \} }{\hat b_i}\right\}.
$$
Then by the ML equation, the monotonicity of $\psi$, and
Lemma~\ref{joskaprop1} of Appendix A, we get
$$
(n-1)M\geq -R_{i_0}=\sum_{j\neq i_0}\psi({\hat a_{i_0}}+{\hat b_j})-(
n-1)\psi({\hat a_{i_0}})\geq
(n-1)[\psi(2{\hat a_{i_0}})-\psi({\hat a_{i_0}})] .
$$
Therefore, ${\hat a_{i_0}}\geq\varepsilon$,
whence ${\hat a_i},{\hat b_i}\geq\varepsilon$ holds
for every $i=1,\dots,n$. 
$\square$

\begin{proposition}\label{joskaprop3} With the solution $\varepsilon$
of $\psi(2x)-\psi(x)=M$ of~(\ref{M}),
we have $f(\varepsilon\1 )\geq \varepsilon\1 $.
\end{proposition}

\noindent \textbf{Proof.}
$$
g_i(\varepsilon\1 )
=\psi^{-1}\left(\psi(2\varepsilon)+\frac{R_i}{n-1}\right)\geq
\psi^{-1}\left(\psi(2\varepsilon)-M\right)=\varepsilon .
$$
 Likewise,
$$
h_i(\varepsilon\1 )
=\psi^{-1}\left(\psi(2\varepsilon)+\frac{C_i}{n-1}\right)\geq
\psi^{-1}\left(\psi(2\varepsilon)-M\right)=\varepsilon .  \quad \square
$$
It is also clear that we have the following.
\begin{proposition}\label{joskaprop4}
If $(\a,\b)\geq (\x,\y )>\0$, then $f(\a,\b)\geq f(\x,\y ).$
\end{proposition}

\begin{theorem}\label{t3}
With $\varepsilon$ satisfying $\psi(2\varepsilon )-\psi(\varepsilon )=M$
 of~(\ref{M}), and  starting at
$\thetav^{(0)} =\varepsilon\1$, the sequence $\thetav^{(it)}$ of the iteration
$\thetav^{(it)}=f(\thetav^{(it-1)})$ for $it\to\infty$ converges at
a geometric rate 
to the unique solution $( {\hat\a},{\hat\b} )$ of the ML equation.
\end{theorem}

\noindent \textbf{Proof.}
From Propositions~\ref{joskaprop2} and~\ref{joskaprop3}
we obtain that the sequence $\thetav^{(it)}$ is coordinate-wise increasing.
Moreover,
it is clear that $(\thetav^{(it)} )$ is bounded from above by
$({\hat\a},{\hat\b})$,
due to Proposition~\ref{joskaprop4}.
Therefore, the convergence of  $\thetav^{(it)}$ follows, and  by the
continuity of $f$, the limit is clearly a fixed point of $f$.
However, in view of Section~\ref{model}, the solution of the ML equation
is a fixed point of $f$, and it cannot be else but the unique solution
$({\hat\a},{\hat\b})$, guaranteed by Theorem~\ref{t1}.
Further, from Theorem~\ref{t2} we get that the rate of convergence is (at least)
geometric. 
$\square$

Therefore, a good starting can be chosen by these considerations. Also note
that at the above $\thetav^{(0)}$ and possibly at its first (finitely many)
iterates,
$f$ is usually not a contraction. It becomes a contraction only when
some iterate $\thetav^{({it}_0 )}$ gets into the neighborhood $K$ of
$\hat \thetav$
of Theorem~\ref{t2}, which is inevitable in view of the convergence of
the sequence $\thetav^{(it)}$. So, Theorem~\ref{t2} is literally applicable
only if we start the iteration at  $\thetav^{({it}_0 )}$. In practice,
however, we do not know the theoretically guaranteed neighborhood $K$.
The practical merit of Theorem~\ref{t3} is just that it offers a
realizable starting.

\section{Applications}\label{application}

First we generated a random directed edge-weighted graph on $n=100$ vertices.
The edge-weight matrix $\W$ had zero diagonal, and the off-diagonal entries
$w_{ij}$'s were independent. Further, for $i\ne j$,
the weight $w_{ij}$ was generated according to
beta-distribution  with parameters $a_i >0$ and $b_j >0$, where $a_i$'s and
$b_j$'s were chosen randomly in the interval [1,5].

Then we estimated the parameters based on $\W$, and plotted the
$a_{i} , {\hat a}_{i}$  $(i=1,\dots ,n)$  and
$b_{j} , {\hat b}_{j}$  $(j=1,\dots ,n)$ pairs.

\begin{figure}[h!]
  \centering
  \includegraphics[width=12cm]{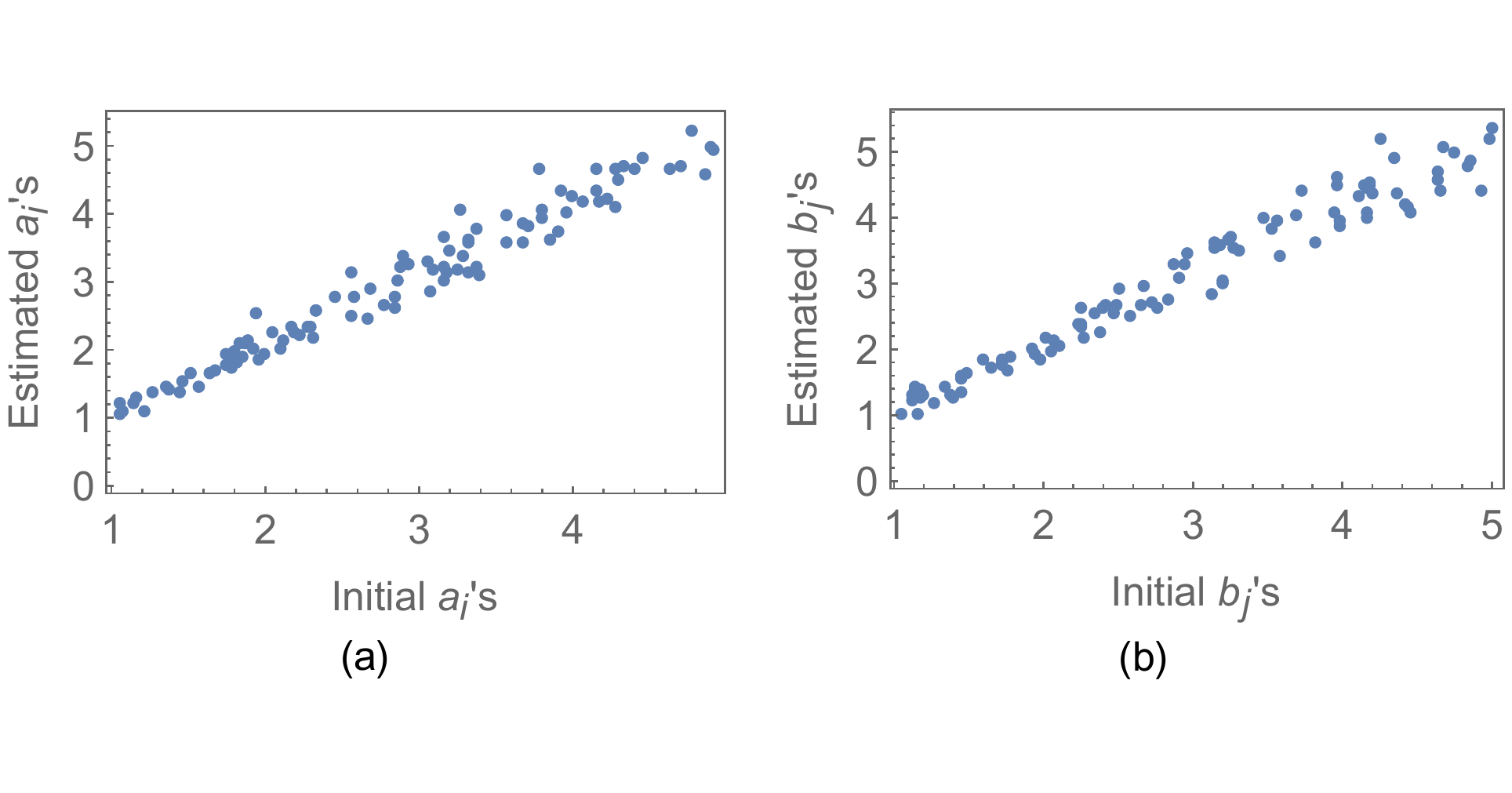}\\
  \caption{Panel (a) shows the original
versus the estimated parameters $a_i$'s with $MSE=0.0628806$, while Panel (b) shows the original
versus the estimated parameters $b_j$'s with $MSE=0.0768382$.} \label{ab}
\end{figure}

Figure~\ref{ab}  shows a good fit between them.

We also applied the
algorithm to migration data between 34 countries. Here $w_{ij}$ is proportional
to the number of people
in thousands who moved from country $i$ to country $j$ (to find jobs)
during the year 2011,
and it is normalized so that be in the interval (0,1).
The estimated parameters are in Table~\ref{migration}.

In this context, $a_i$'s are related to the  emigration and
and $b_i$'s to the counter-immigration potentials.
When $a_i$ is large, country $i$ has a relatively large potential
for emigration. On the contrary, when
$b_i$ is large, country $i$ tends to have a relatively large
resistance against immigration.

\begin{table}[h!]
\scriptsize
\begin{tabular}{c l c r @{.} l@{\hskip 0.35in} c l c r @{.} l }
  \hline
  i & Country & $a_i$  & \multicolumn{2}{@{\hskip -0.3in}c}{$b_i$} & i & Country & $a_i$ & \multicolumn{2}{c}{$b_i$}  \\\hline
  1 & Australia & 0.26931 & 1475&75242 & 18 & Japan & 0.23211 & 9926&91644 \\
  2 & Austria & 0.27403 & 632&81653 & 19 & Korea & 0.22310 & 4199&25005 \\
  3 & Belgium & 0.33380 & 46&01197 & 20 & Luxembourg & 0.17543 & 107&91399 \\
  4 & Canada & 0.27383 & 2363&23435 & 21 & Mexico & 0.26706 & 4655&95370 \\
  5 & Chile & 0.21236 & 28940&59777 & 22 & Netherlands & 0.37754 & 39&52320 \\
  6 & Czech Rep.& 0.31188 & 470&28651 & 23 & New Zealand & 0.20542 & 2568&00582 \\
  7 & Denmark & 0.26514 & 847&34887 & 24 & Norway & 0.22646 & 519&12451 \\
  8 & Estonia & 0.23235 & 25602&33371 & 25 & Poland & 0.62846 & 1106&55946 \\
  9 & Finland & 0.29357 & 1100&00568 & 26 & Portugal & 0.31011 & 1606&59979 \\
  10 & France & 0.52721 & 37&92122 & 27 & Slovak Rep. & 0.27871 & 42451&19093 \\
  11 & Germany & 0.62020 & 1&64064 & 28 & Slovenia & 0.19720 & 6824&54028 \\
  12 & Greece & 0.29708 & 6319&19184 & 29 & Spain & 0.39732 & 182&47160 \\
  13 & Hungary & 0.31443 & 32750&88310 & 30 & Sweden & 0.39627 & 57&34509 \\
  14 & Iceland & 0.18051 & 2950&72653 & 31 & Switzerland & 0.33611 & 4524&67821 \\
  15 & Ireland & 0.27555 & 364&52781 & 32 & Turkey & 0.25900 & 146175&82805 \\
  16 & Israel & 0.25854 & 1926&04551 & 33 & United Kingdom & 0.49301 & 48&61626 \\
  17 & Italy & 0.50522 & 135&14076 & 34 & United States & 0.38019 & 2433&78269 \\
  \end{tabular}
\caption{Estimated parameters for migration data, 2011}
\label{migration}
\end{table}

It should be noted again that edge-weighted graphs of this type very
frequently model real-world directed networks.

\section*{Appendix}\label{App}

\subsection*{A.  Properties of the digamma function}

Though, we do not use it explicitly, the following approximation of the
digamma function  $\psi (x)=\frac{\partial \ln \Gamma (x)}{\partial x} =
\frac{\Gamma' (x)}{\Gamma (x)}$ $(x>0$)
is interesting for its own right.

\begin{lemma}\label{Taylor}
$\psi (x) = \ln (x-\frac12 ) + {\cal O} \left( \frac1{x^2 } \right)$ for $x>1$.
\end{lemma}
The statement of the lemma easily follows by Taylor expansion.

\begin{lemma}\label{psi1}
$\frac1{\psi'(x+y)}>\frac1{\psi'(x)}+\frac1{\psi'(y)}$ for $x,y>0$.
\end{lemma}

\noindent
\textbf{Proof.}
First we prove that the function
$u(x)=\frac1{\psi'(x)}$, $x\in (0,\infty)$ is strictly
convex. Indeed,  one can easily see that
$$
 u''(x)=\frac{-\psi'''(x)[\psi'(x)]^2+2[\psi''(x)]^2\psi'(x)}{[\psi'(x)]^4} ,
$$
and this is positive due to $\psi' (x)>0$ and the fact that
$$
 \frac{[\psi''(x)]^2}{\psi'''(x)\psi'(x)}>\frac12 .
$$
Latter one is a particular case
of Corollary 2.3 in~\cite{alzer}.

Now, in view of $\lim_{x\to0}\psi'(x)=\infty$, we can extend $u$ continuously
to 0 by setting
$u(0)=0.$ Then $u$ is still strictly convex,
and therefore, for every $x,y>0$ we have
$u(x)=u(\frac{y}{x+y}\cdot 0+\frac{x}{x+y}\cdot(x+y))<
\frac{y}{x+y}u(0)+\frac{x}{x+y}u(x+y)$.
Consequently,
\begin{equation}\label{convineq1}
u(x)< \frac{x}{x+y} u(x+y),
\end{equation}
and likewise,
\begin{equation}\label{convineq2}
u(y)< \frac{y}{x+y}u(x+y).
\end{equation}
Adding (\ref{convineq1}) and (\ref{convineq2}) together,
we get the statement of the lemma. $\square$

\begin{lemma}\label{joskaprop1}
The function $\psi(2x)-\psi(x)$, $x\in (0,\infty)$
is decreasing and its range is $(\ln 2,\infty)$.
\end{lemma}

\noindent \textbf{Proof.}
It is easily seen by the identity
$\psi(2x)=\frac12\psi(x)+\frac12\psi\left(x+\frac12\right)+\ln 2$
which can be found in~\cite{handbook}. $\square$

In the last lemma we collect some limiting properties of the digamma function
and its derivative, see, e.g.,~\cite{handbook,alzer,Bernardo} for details.

\begin{lemma}\label{ditul}
The digamma function $\psi$ is a strictly concave, smooth function on
$(0,\infty)$ that satisfies the following limit relations:
$$
\begin{aligned}
 &\lim_{x\to 0+} \psi (x) =-\infty, \quad
 \lim_{x\to \infty } \psi (x) =\infty, \quad
 \lim_{x\to \infty } \psi' (x) =0, \\
 &\lim_{x\to \infty } (\psi (x) -\ln x ) =0, \quad \quad
 \lim_{x\to 0+} (\psi (2x) -\psi (x)) =\infty .
\end{aligned}
$$
\end{lemma}

\subsection*{B. Considerations on the boundary of the mean parameter space}

In Section~\ref{model}, we saw that the mean parameter space ${\cal M}$
consists of $2n$-tuples
$(A_1 ,\dots A_n ,B_1 ,\dots ,B_n)$ obtained from the parameters
$(\a ,\b ) = (a_1 ,\dots a_n ,b_1 ,\dots ,b_n)$ of the underlying
beta-distributions by Equations~(\ref{mean}).

Denoting by $L (\a ,\b ) =(A_1 (\a, \b ) , \dots ,A_n (\a ,\b ),
B_1 (\a ,\b ),\dots , B_n (\a ,\b))$ this dependence, i.e., the
$\Theta \to {\cal M}$ (one-to-one) mapping, a boundary
point $\bar L =({\bar A}_1 ,\dots ,{\bar A}_n ,{\bar B}_1 ,\dots ,{\bar B}_n)$
of $\cal M$ can be obtained as
${\bar L} =\lim_{k\to \infty} L (\a^k ,\b^k )$, where
$\a^k = (a^k_{1} ,\dots ,a^k_{n} )$, $\b^k = (b^k_{1} ,\dots ,b^k_{n} )$, and
$$
\begin{aligned}
 \lim_{k\to \infty } a^k_{i} &= {\bar a}_i \in [0,\infty], \quad i=1,\dots n; \\
 \lim_{k\to \infty } b^k_{j} &= {\bar b}_j \in [0,\infty], \quad j=1,\dots n .
\end{aligned}
$$
In view of Lemma~\ref{ditul}, only the
${\bar a}_i,{\bar b}_j =\infty$ cases have relevance. 
The sequence $(\a^k ,\b^k )$ can be chosen such that
\begin{equation}\label{tart}
 \lim_{k\to\infty } \frac{a^k_{i}}{b^k_{j}} =x_{ij}  \quad \textrm{with} \quad
 0<x_{ij} < \infty , \quad  \textrm{for} \quad  i\ne j.
\end{equation}
Then, using~\ref{mean}),
$$
\begin{aligned}
 {\bar A}_i &= -\lim_{k\to\infty} \sum_{j\ne i} \left[ \psi ({a}^k_{i}+
  {b}^k_{j})- \psi ({a}^k_{i} ) \right]
 =\lim_{k\to\infty} \sum_{j\ne i} \left[ \ln ({a}^k_{i} ) -
  \ln ({a}^k_{i}  +{b}^k_{j}) \right] =
 \sum_{j\ne i} \ln \frac{x_{ij}}{1+x_{ij}},  \quad i=1,\dots ,n; \\
 {\bar B}_j &= -\lim_{k\to\infty} \sum_{i\ne j} \left[ \psi ({a}^k_{i}+
   {b}^k_{j}) - \psi ({b}^k_{j} ) \right]
 =\lim_{k\to\infty} \sum_{i\ne j} \left[ \ln ({b}^k_{j} ) -
  \ln ({a}^k_{i}  +{b}^k_{j}) \right] = \sum_{i\ne j} \ln \frac{1}{1+x_{ij}} ,
 \quad j=1,\dots ,n .
\end{aligned}
$$
These equations show that the boundary point $\bar L$ of $\cal M$
contains -- in its coordinates --
the row- and column-sums of the matrices $\U (\W )$ ad $\V (\W )$ respectively
(see Section~\ref{model}), where the general off-diagonal entry of the
 $n\times n$ edge-weight matrix $\W$ is $\frac{x_{ij}}{1+x_{ij}}$.

Observe that $2n-1$ $x_{ij}$'s can be chosen free, and all the others
are obtainable from them. To see this, consider the complete bipartite graph
on vertex classes $(a_1 ,\dots ,a_n)$ and $(b_1 ,\dots ,b_n)$,
where to the edge connecting $a_i$ and $b_j$ we assign $\frac{a_i }{b_j}$.
Choose a minimal spanning tree of this graph (it contains $2n-1$ edges), and
consider the sequence of $(\a^k ,\b^k )$'s satisfying condition~(\ref{tart}).
Then, as $k\to \infty$,
the  $x_{ij}$'s of the edges not included in the spanning tree
can be obtained from the $x_{ij}$'s  of the $2n-1$
edges included in  the spanning tree. Therefore the row- and
column-sums of the  edge-weight matrix $\W$ of entries
$\frac{x_{ij}}{1+x_{ij}}$ $(i\ne j)$
are on a $(2n-1)$-dimensional manifold
in $\R_{-}^{2n}$, so they are on the boundary of the convex hull $\cal T$
of the possible sufficient statistics $(\RR ,\CC )$.
However,
this boundary has zero Lebesgue measure, and so, zero probability
with respect to the underlying absolutely continuous distribution.

\subsection*{C.  Proof of Theorem~\ref{t2}}

It suffices to prove that some induced matrix norm of the matrix of the
first  derivatives $\J$ of $f$ at $\hat \thetav$ is strictly less than 1.
We prove this for the $L_1$-norm. From~(\ref{fixes}) we obtain that
\begin{equation}\label{parcszuk}
\frac{\partial g_i}{\partial a_i} ({\hat \a }, {\hat b} ) =
\frac {\displaystyle \frac1{n-1}\sum_{j\neq
i}\psi'({\hat a}_i+{\hat b}_j)}
 {\displaystyle\psi'\left[ \psi^{-1}\left(\frac1{n-1}\sum_{j\neq i}
\psi({\hat a}_i+{\hat b}_j)+\frac{R_i}{n-1}\right)\right] }.
\end{equation}
From~(\ref{le}) we have $\frac1{n-1}\sum_{j\neq i}
\psi({\hat a}_i+{\hat b}_j)+\frac{R_i}{n-1}=\psi({\hat a}_i)$.
Substituting it into~(\ref{parcszuk}), we get
$$
\frac{\partial g_i}{\partial a_j}({\hat \a},{\hat b})= \left\{
\begin{array}{ll}
\frac {\displaystyle \sum_{s\neq i} \frac1{n-1}\psi'({\hat a}_i+{\hat b}_s)}
{\displaystyle\psi'({\hat a}_i)} & \mbox{\ if\ } j=i\\
0 & \mbox{\ if\ } j\neq i .
\end{array}\right.
$$
Likewise,
$$
\frac{\partial g_i}{\partial b_j}({\hat \a},{\hat b})= \left\{
\begin{array}{ll}
0 & \mbox{\ if\ } j=i\\ \frac {\displaystyle \frac1{n-1}
 \psi'({\hat a}_i+{\hat b}_j)} {\displaystyle\psi'({\hat a}_i)}
& \mbox{\ if\ } j\neq i .
\end{array}\right.
$$
Further,
$$
\frac{\partial h_i}{\partial a_j}({\hat \a},{\hat b})= \left\{
\begin{array}{ll}
0 & \mbox{\ if\ } j=i\\
\frac {\displaystyle \frac1{n-1}\psi'({\hat a}_j+{\hat b}_i)}
{\displaystyle\psi'({\hat b}_i)} & \mbox{\ if\ } j\neq i
\end{array}\right.
$$
and
$$
\frac{\partial h_i}{\partial b_j} ( {\hat \a},{\hat \b} )= \left\{
\begin{array}{ll}
\frac {\displaystyle \sum_{s\neq i} \frac1{n-1}\psi'({\hat a}_s+{\hat b}_i)}
{\displaystyle\psi'({\hat b}_i)} & \mbox{\ if\ } j=i\\
0 & \mbox{\ if\ } j\neq i .
\end{array}\right.
$$
Observe that $J ({\hat \a} ,{\hat b} )$ has nonnegative entries. Therefore,
its $L_1$-norm is the maximum of its column-sums.
The $j$th column-sum of $J ({\hat \a} ,{\hat b} )$ is equal to
\begin{equation}\label{colsum}
\frac{\displaystyle \sum_{s\neq j} \frac1{n-1}\psi'({\hat a}_j+{\hat b}_s)}
{\displaystyle\psi'({\hat a}_j)}
+ \sum_{s\neq j}\frac {\displaystyle \frac1{n-1}\psi'({\hat a}_j+{\hat b}_s)}
{\displaystyle\psi'({\hat b}_s)}
=\frac1{n-1}\sum_{s\neq j}\displaystyle \psi'({\hat a}_j+{\hat b}_s)\left(
\displaystyle\frac1{\psi'({\hat a}_j)}+\frac1 {\displaystyle \psi'({\hat b}_s)}
 \right)
\end{equation}
for $j=1,\dots ,n$; and likewise, the $(n+j)$th column-sum of
$J ({\hat \a} ,{\hat b} )$ is
\begin{equation}\label{colsum1}
\frac1{n-1}\sum_{s\neq j}\displaystyle \psi'({\hat a}_s+{\hat b}_j)\left(
\displaystyle\frac1{\psi'({\hat a}_s)}+\frac1 {\displaystyle \psi'({\hat b}_j)}
 \right)
\end{equation}
for $j=1,\dots ,n$.
As (\ref{colsum}) and (\ref{colsum1}) are of similar appearance, it suffices
to prove that the right hand side of  (\ref{colsum}) is less than 1.
But $\displaystyle \psi'({\hat a}_j+{\hat b}_s)\left(
\displaystyle\frac1{\psi'({\hat a}_j)}+\frac1 {\displaystyle
\psi'({\hat b}_s)}\right)<1$ holds by Lemma~\ref{psi1},
and we have $n-1$ terms in the summation.

Since $f: \R^{2n} \to \R^{2n}$ is continuously differentiable in a
neighborhood of ${\hat \thetav }=({\hat \a} ,{\hat \b})$, Theorem 3
of~\cite{Grasmair} implies that there is a closed neighborhood $K$ of
$\hat \thetav$ such that $f$ is a contraction on $K$.
In particular, the fixed point iteration $f (\thetav^{(it -1)} )=
\thetav^{(it )}$ $(it\to \infty )$
converges for every $\thetav^{(0 )} \in K$ to $\hat \thetav$,
which is the unique solution of~(\ref{fixes}). $\square$

\end{document}